\documentclass[psamsfonts]{amsart}
\usepackage{amsfonts,amssymb,amscd,amsmath,bbding,enumerate,verbatim,xypic, mathdots,mathtools}
\usepackage{epsfig}
\usepackage[colorlinks=true
]{hyperref}

\def\lto{{\longrightarrow}}
\def\into{{\hookrightarrow}}
\def\xto{\xrightarrow}

\newcommand{\calh}{{\mathcal H}}

\newcommand{\call}{{\mathcal L}}
\newcommand{\calm}{{\mathcal M}}

\newcommand{\calo}{{\mathcal O}}
\newcommand{\calp}{{\mathcal P}}

\newcommand{\calt}{{\mathcal T}}

\newcommand{\calox}{{\calo_X}}

\def\calend{{\mathcal E}nd\,}

\def\calhom{{\mathcal H}om}

\newcommand{\bbbr}{{\mathbb R}}

\newcommand{\CC}{{\mathbb C}}
\newcommand{\DD}{{\mathbb D}}
\newcommand{\EE}{{\mathbb E}}

\newcommand{\LL}{{\mathbb L}}

\newcommand{\PP}{{\mathbb P}}

\newcommand{\RR}{{\mathbb R}}

\newcommand{\VV}{{\mathbb V}}
\newcommand{\ZZ}{{\mathbb Z}}

\newcommand{\bdot}{\bullet}

\DeclareMathOperator{\ev}{ev}
\DeclareMathOperator{\End}{ End}

\DeclareMathOperator{\Ext}{Ext}

\DeclareMathOperator{\Hoch}{HH}

\DeclareMathOperator{\Hom}{Hom}

\DeclareMathOperator{\id}{id}

\DeclareMathOperator{\Ker}{Ker}

\DeclareMathOperator{\op}{op}

\DeclareMathOperator{\QCoh}{\mathsf {QCoh}}
\DeclareMathOperator{\rad}{rad}

\DeclareMathOperator{\RHom}{\bbbr\Hom}

\DeclareMathOperator{\Spec}{Spec}

\DeclareMathOperator{\Sym}{{\mathbb S}ym}

\DeclareMathOperator{\Tor}{Tor}

\DeclareMathOperator{\Coh}{\mathsf{Coh}}
\DeclareMathOperator{\fmod}{\ensuremath{ \mathsf{mod}}}

\DeclareMathOperator{\Loc}{\ensuremath{ \mathsf{Loc}}}

\theoremstyle{definition}
\newtheorem{defn}{Definition}[section]

\newtheorem{question}[defn]{Question}

\newtheorem{remark}[defn]{Remark}

\newtheorem{sit}[defn]{}

\newtheorem{example}[defn]{Example}

\theoremstyle{plain}

\newtheorem{proposition}[defn]{Proposition}

\newtheorem{theorem}[defn]{Theorem}

\newtheorem{lemma}[defn]{Lemma}

\newtheorem{cor}[defn]{Corollary}

\theoremstyle{remark}

\begin{document}
\title[Hochschild (Co-)homology and Tilting] 
{Hochschild (Co-)Homology  of Schemes with Tilting Object}


\author[R.-O.~Buchweitz]{Ragnar-Olaf~Buchweitz}
\address{Department of Computer and Mathematical Sciences,
University of Toronto Scarborough, Toronto, Ontario,
Canada M1C 1A4} \email{ragnar@utsc.utoronto.ca}

\author[L.~Hille]{Lutz Hille}
\address{Mathematisches Institut
der Universit\"at M\"unster,
Einsteinstra\ss e 62,
48149 M\"unster,
Germany} 
\email{lutzhille@uni-muenster.de}

\thanks{The first author gratefully acknowledges partial support through 
NSERC grant 3-642-114-80, while the second author thanks
SFB 478 "Geometrische Strukturen in der Mathematik" for its support.}

\date{\today} 

\begin{abstract}
Given a $k$--scheme $X$ that admits a tilting object $T$, we prove that the Hochschild 
(co-)homology of $X$ is isomorphic to that of $A=\End_{X}(T)$. We treat more generally 
the relative case when $X$ is flat over an affine scheme $Y=\Spec R$ and the tilting object 
satisfies an appropriate Tor-independence condition over $R$.
Among applications, Hochschild homology of $X$ over $Y$ is seen to vanish in negative 
degrees, smoothness of $X$ over $Y$ is shown to be equivalent to that of $A$ over $R$, 
and for $X$ a smooth projective scheme we obtain that Hochschild homology is concentrated 
in degree zero. Using the Hodge decomposition \cite{BFl2} of Hochschild homology 
in characteristic zero, for $X$ smooth over $Y$ the Hodge groups 
$H^{q}(X,\Omega_{X/Y}^{p})$ vanish for $p < q$, while in the absolute case
they even vanish for $p\neq q$.

We illustrate the results for crepant resolutions of quotient singularities, in particular for 
the total space of the canonical bundle on projective space.
\end{abstract}

\maketitle
{\footnotesize\tableofcontents}

\section*{Introduction}
If $A$ and $B$ are derived equivalent algebras over a field $k$, then their Hochschild theories 
are isomorphic; see e.g. \cite{Kel3}. If $X$ and $Y$ are smooth complex projective varieties that 
are derived equivalent through a Fourier--Mukai transformation, their Hochschild theories agree 
as well; see e.g. \cite{Cal}.

Here we consider the case when a scheme $X$ is derived equivalent to an algebra $A$
through a {\em classical tilting object\/} in its derived category of quasi-coherent sheaves
and establish that the Hochschild theories of $X$ and that of $A$ are naturally isomorphic too. 

More precisely, for $k$ a field, we assume that $X$ is a $k$--scheme projective over an affine 
scheme $Z=\Spec K$ of finite type over $k$ and that $T$ is a tilting object on $X$. 
In that case, the key theorem of geometric tilting theory, recalled in Theorem \ref{projtilt} below, 
yields through $\RR\Hom_{X}(T,?)$ an exact equivalence from the derived category of 
quasicoherent $\calox$--modules to the derived category of (right) $A=\End_{X}(T)$--modules. 

Our main results, Corollary \ref{Hochschild CoInvariance} and Theorem \ref{invariance complex},
are now that if such a scheme $X$ is flat over an affine $k$--scheme ${Y=\Spec R}$ so that the endomorphism ring $A$ of the tilting object is as well flat over $R$, then the 
Hochschild (co-)homology of $X$ over $Y$  becomes naturally isomorphic to the Hochschild (co-)homology of $A$ over $R$. Here, in accordance with the results in \cite{Swan, BFl1}, 
Hochschild theory for a flat morphism $X\to Y$ is defined to be the hyper-\-(co-)homology 
theory attached to $\Delta_{*}\calox$, where $\Delta:X\to X\times_{Y}X$ is the diagonal embedding,
while Hochschild theory of $A$ over $R$ is understood as the (co-)homology theory attached to 
$A$ as a (right) module over the enveloping algebra $A^{\op}\otimes_{R}A$, which coincides with
Hochschild's original definition when $A$ is projective as $R$--module, and in the flat case
is the specialization of Quillen's approach in \cite{Qu}.

After reviewing briefly the theory and scope of tilting objects in algebra and geometry 
in Section 1, we investigate appropriate Tor-independence properties of tilting objects, such 
as flatness, in Section 2, to establish the invariance of Hochschild cohomology in Section 3, 
that of Hochschild homology in Section 4. Section 5 uses the Hodge decomposition for Hochschild homology in characteristic zero to obtain vanishing of various Hodge groups in the context of tilting and
Section 6 deals with the example of crepant resolutions of some quotient singularities.

As applications we can strengthen some earlier results. For example; see 
Corollary \ref{homology vanishing}; existence of a tilting object forces $H^{i}(X,\calox)=0$ for 
$i\neq 0$ independent of the characteristic of $k$, and for $X$ flat over $Y$, the {\em negative\/} Hochschild homology, $\Hoch_{i}(X/Y)$ for $i<0$, vanishes.

For smooth projective schemes, we even obtain $\Hoch_{i}(X)=0$ for $i\neq 0$; see 
Theorem \ref{vanishing absolute}. Further, while smoothness of $X$ was known to imply finite 
global dimension for $A$, here we show in Corollary \ref{smoothness} that for $X$ and its 
tilting object $T$ flat over $Y$, smoothness of $X$ over $Y$ is indeed equivalent to that of 
$A$ over $R$, where smoothness in the noncommutative algebraic context is interpreted in 
the sense of van den Bergh \cite{vdB}. 

Employing the Hodge decomposition of Hochschild cohomology in characteristic zero 
\cite{Swan, BFl2}, it follows for a smooth morphism $X\to Y$ over the complex numbers 
that  the Hodge groups $H^{q}(X,\Omega_{X/Y}^{p})$ vanish for $p < q$, while in the absolute 
case even  $H^{q}(X,\Omega_{X}^{p})=0$ for $p\neq q$.

In the final section, we illustrate the results for the total space of the canonical bundle over projective 
space, a crepant resolution of a quotient singularity for an algebraically closed field of 
characteristic zero.

We point out that in particular examples, at least in the absolute case of complex smooth projective schemes, invariance of Hochschild (co-)homology or partial consequences thereof
have already been known or alluded to by some authors; see, for example, 
\cite[1.8.2]{BMa}, \cite[2.2]{CGW} or the introduction of \cite{SWi}.

To end this introduction, we comment on the broader picture. As detailed in 
\cite[Thm. 7.5]{HvdB} on any separated scheme $X$ there exists a perfect complex 
$E$ such that the derived category of quasicoherent $\calox$--modules is equivalent 
to the derived category of the differential graded algebra (DG algebra) $A=\RR\Hom_{X}(E,E)$. 
To a category such as the latter, To$\ddot{\text{e}}$n \cite[8.1]{Toen} assigns a Hochschild theory that is 
essentially intrinsic in his context, thus, should specialize both to the geometric incarnation of the Hochschild theory of $X$ on the one hand and to the algebraic realization of the 
Hochschild theory of $A$ on the other when $E$ is a classical tilting object.

In the same vein, the flatness or Tor-independence conditions imposed here should be 
avoidable if one employs the theory of DG algebras and resolvents for morphisms of 
schemes or analytic spaces as was done in \cite{BFl1}. 
We decided, however, to present the ``classical'' version with a direct proof that
avoids the formidable technical apparatus certainly necessary for the ultimate treatment 
of invariance of Hochschild theory under a larger class of exact equivalences.

\section{Classical Tilting Objects}
Let $k$ be a field and $\calt$ a triangulated $k$--linear category. 
With $[i]$, as usual, the $i^{th}$ iteration of the given (translation) auto-equivalence 
on $\calt$, we set $\Ext^{i}_{\calt}(M,N)=\Hom_{\calt}(M,N[i])$.
Recall that a subcategory of $\calt$ is {\em thick\/}, if it is closed under translations, 
exact triangles and direct summands, and {\em localizing\/} if it is further closed under 
all {\em small\/}, that is, set-indexed, direct sums that exist in $\calt$. If $U$ is any object or (full) subcategory of $\calt$, we denote $\Loc U$ the smallest localizing subcategory containing 
$U$ in $\calt$.

Concerning functors between triangulated $k$--linear categories we only allow those that
are {\em $k$--linear\/} and {\em exact\/}.

\subsection*{Tilting Objects in Triangulated Categories}
We first introduce the concept of a tilting object in a ``large'' triangulated category. The reason 
for this is that there the generating condition has a formulation that is easier to verify even when 
one is ultimately only interested in triangulated categories whose size is bounded in some way.
The relevant definition is the following.

\begin{defn}
\label{defn:tilt}
Let $\calt$ be a triangulated category that is closed under small direct sums. 
An object $T$ in $\calt$ is {\em tilting\/} if it is a {\em compact generator with
only trivial self-extensions\/},
that is to say
\begin{enumerate}[\quad\rm (1)]
\item ({\em Compactness\/}) The functor $\Hom_{\calt}(T,?)$ commutes with small direct sums;
\item
\label{defn:tilt:gen} 
({\em Generating\/}) We offer two versions that are equivalent in the presence of (1):
\begin{enumerate}[\rm (a)]
\item If $N$ in $\calt$ satisfies $\Ext^{i}_{\calt}(T,N)=0$ 
for each $i\in\ZZ$, then $N=0$;
\item 
\label{defn:tilt:genb} 
The smallest localizing subcategory that contains $T$ is $\calt$.
\end{enumerate}
\item ({\em  Only trivial self-extensions\/}) $\Ext^{i}_{\calt}(T,T)=0$ for $i\neq 0$.
\end{enumerate}
\end{defn}

\begin{remark}
As to the equivalence of (a) and (b) above, note that for any object $N$ in $\calt$ 
the full subcategory $^{\perp}N$ consisting of those objects $X$ from $\calt$ with 
$\Ext^{i}_{\calt}(X,N)=0$  for each $i\in\ZZ$ is localizing. Therefore, (b) $\Rightarrow$ (a), 
as $T$ in $^{\perp}N$ then implies $\Loc T=\calt\subseteq {^{\perp}N}$, which means 
that the identity morphism on $N$ is zero. 

The converse requires $T$ to be compact, thus (1). Namely, if $T$ in $\calt$ satisfies 
(1) and (2a), then the set of objects $\bigcup_{n\in\ZZ}T[n]$ {\em compactly generates\/} 
$\calt$ in the sense of definition \cite[1.7]{Nee}, and Theorem 2.1.2. in that reference gives (a) 
$\Rightarrow$ (b).
\end{remark}

\begin{remark}
Tilting objects as just defined are nowadays sometimes called ``classical'' to distinguish them
from the more general notion where $\calt$ carries a DG enhancement, and rather than requiring 
no self-extensions one considers the full DG algebra of endomorphism in the enhancement.
In that case, the target category becomes that of DG modules over the DG algebra.
See \cite{HvdB} for further details.
\end{remark}

\begin{sit}
\label{commentsgen}
Despite appearances, (b) is often easier to verify than (a). If we can show that a known (set of) 
compact generator(s) is contained in $\Loc T$, then $T$ is already generating by (b). For 
example, in practice it is often known beforehand that the triangulated category $\calt$ in 
question is {\em compactly generated\/}, in the sense that there is a set of 
compact objects $G$ that satisfies (a) or, equivalently, (b). To test then that a given
object $T$ is compact and generates, it suffices to check
\begin{enumerate}[\qquad\rm (c)]
\item 
\label{genc} The smallest thick subcategory in $\calt$ that contains $T$ equals the
subcategory $\calt^{c}$ of all compact objects.
\end{enumerate}

Another way to employ the equivalence of the generating conditions is as follows. 
Assume $T$ in $\calt$ satisfies (1) and (2) and let $L:\calt\to\calt'$ be an exact functor into 
a triangulated category also closed under small direct sums,  and assume that $L$ commutes 
with small direct sums, for example, if $L$ admits a right adjoint. It then follows that in 
$\calt'$ the full subcategory $L(\calt)$ is contained in $\Loc L(T)$. Thus, if $L(\calt)$ contains 
a generating set of compact objects for $\calt'$, then $L(T)$ satisfies (1) and (2) in $\calt'$ as 
soon as it is again compact, and only (3) remains to be verified to establish $L(T)$ as 
a tilting object in $\calt'$.
\end{sit}

\begin{example}
\label{algtilt}
(See \cite{Kel, Kel3} for a more general account and further references.)
For a ring $B$, denote $D(B)$ the full derived category of right $B$--modules.

Assume it is known that the triangulated category $\calt$ is equivalent to the derived 
category of some {$k$--algebra\/} $B$.  In that case $\calt$ must contain a tilting object $T$, 
as any ring, when considered as a module over itself, is a tilting object in its own derived category. 

By Rickard's fundamental result, for any tilting object $T$, the category $\calt$ is 
equivalent to $D(\End_{\calt}(T))$. To be more precise, assume there is an
equivalence $F\colon\calt \to D(B)$. As any equivalence, the functor $F$ will preserve compactness, 
whence $F(T)$ is a {\em perfect\/} complex of $B$--modules, as those complexes are 
precisely the compact objects in $D(B)$; see e.g. \cite[Prop.9.6]{Chr} for a proof. 
As the generating property and lack of self-extensions are as well preserved by the 
equivalence $F$, Rickard tells us that $?\otimes^{\LL}_{B}F(T)$ provides an equivalence 
from $D(B)$ onto $D(A')$, where $A'=\End_{B}(F(T))$. Thus, $F(?)\otimes^{\LL}_{B}F(T)$
provides an equivalence from $\calt$ to $D(A')$. Finally note that $F$ induces an isomorphism 
of algebras $A=\End_{\calt}(T)\cong A'=\End_{B}(F(T))$, whence in summary, 
$\calt \cong D(\End_{\calt}(T))$ as claimed.

It follows that tilting objects {\em detect all $k$--algebras\/} $B$ that satisfy $D(B)\cong \calt$ as triangulated categories in the sense that the assignment $T\mapsto \End_{\calt}(T)$ is a surjection 
onto the isomorphism classes of those rings $B$. It would be interesting to understand the fibres 
of this assignment. Note, for example, that the {\em Picard group\/} of $\calt$, that is, the group of 
all auto-equivalences of $\calt$, operates on those fibres. 

Derived equivalent algebras have isomorphic Hochschild (co-)homology, which thus becomes an invariant of the derived category $\calt$, retrievable 
from the endomorphism algebra of any tilting object.
\end{example}

\subsection*{Tilting Objects in Geometry}
Here we are mainly interested in the situation, where $\calt$ is a  ``{\em geometric\/}'' 
triangulated category, in that  $\calt \cong D(X) = D(\QCoh(X ))$, the derived 
category of quasi-coherent sheaves on some scheme $X$ over $k$.
For any noetherian quasi-projective scheme, the triangulated category $D(X)$ 
is closed under small direct sums; see e.g. \cite[Example 1.3.]{Nee}  and the 
references therein. Moreover, that category is compactly generated, and the compact objects
are exactly the perfect complexes, as soon as $X$ is quasi-compact and separated; 
see \cite[Prop.~2.5.]{Nee}.

\begin{sit}
\label{class}
The category of schemes we will consider consists of those 
$k$--schemes $X$ such that the structure morphism $X\to \Spec k$ can be factored as 
$X\xto{p}Z\xto{q}\Spec k$ with $p$ projective and $Z$ an affine scheme of finite type 
over $k$. The morphisms between such schemes are the morphisms over $\Spec k$. 

For such a scheme $X$, the triangulated category $D(X)$ is thus in particular
$k$--linear, closed under small direct sums, and compactly generated.
\end{sit}

\begin{remark}
\label{virtproj}
Let us point out that the affine scheme $Z$ appearing above plays only an 
auxiliary role. Indeed, assume $Y$ is any affine scheme over $k$ of finite type and suppose 
the structure morphism $X\to \Spec k$ factors through a morphism $f:X\to Y$ as well 
as through a projective morphism $p:X\to Z$ as before. Then $Z'=Y\times_{k}Z$ 
is again affine of finite type over $k$ and the induced morphism $X\to Z'$ factors $p$, thus, 
is in turn projective; see \cite[Prop.5.5.5]{EGAII}.

This flexibility implies, for example, that the category of schemes under consideration
is closed under fibre  products over affine schemes of finite type, so that for $X, X'$ in that 
category and $f:X\to Y, f':X'\to Y$ morphisms to an affine scheme $Y$ of finite type over $k$, 
the fibre product $X\times_{Y}X'$ belongs again to the category.

In particular, the category is closed under base change by morphisms  $Y'\to Y$ of affine 
schemes of finite type over $k$, that is, with $f:X\to Y$ a morphism, $X'=X\times_{Y}Y'$ is 
again in that category if $X$ is.
\end{remark}

The reason to restrict ourselves to the category of schemes in \ref{class} is the following 
structural result that has its origin in Beilinson's seminal paper \cite{Bei} and was 
developed further through \cite{Baer}, \cite[Thm.6.2]{Bon1} and \cite{BvdB}. 
The form given here is \cite[Thm.7.6]{HvdB}. To abbreviate, we call $T$ from $D(X)$
a {\em tilting object on $X$\/} if it is one for that triangulated category.

\begin{theorem}
\label{projtilt}
If $X$ as in {\em \ref{class}\/} admits a tilting object $T$, then with $A=\End_{X}(T)$ 
the following hold:
\begin{enumerate}[\rm\quad(1)]
\item 
The functor $T_{*}=\RHom_{\calox} (T,\ )$ induces an equivalence from $D(X)$ to 
$D(A)$. Its left adjoint $T^{*}=(\ )\otimes^{\LL}_{A}T$ provides the inverse equivalence.
\item
The equivalence $T_{*}$ carries $D^{b}(\Coh(X))$, the bounded derived category of coherent 
$\calox$--modules, to $D^{b}(\fmod A)$, the bounded derived category of finitely 
generated right $A$--modules.
\item
\label{fgd}
If $X$ is smooth then $A$ has finite global dimension.
\end{enumerate}
\end{theorem}

\begin{remark}
\label{finalg}
If $X$ maps to some affine $k$--scheme $Y=\Spec R$, then the equivalence $T_{*}$ 
is $R$--linear, thus, $A$ is naturally an $R$--algebra. If $X\xto{p}Z=\Spec K \xto{q}\Spec k$ 
is a factorization with $p$ projective, as is supposed to exist, then, in view of 
\cite[Thm. 2.4.1.(i)]{EGAIII}, the ring $A$ is a {\em finite\/} $K$--algebra. In particular, as a 
ring, $A$ is noetherian on either side and (module--)finite over its centre.
\end{remark}

We now turn to some basic examples.

\subsection*{Tilting in the Absolute Case}
\begin{sit}
If $X$ is already projective over the field $k$, and $T$ a tilting object on it, 
then $A=\End_{X}(T)$ is a finite-dimensional $k$--algebra
and so its Grothendieck group $K_{0}(A)$ of finitely generated modules is free abelian 
of finite rank.  In view of the equivalence, this is then also isomorphic to the Grothendieck 
group $K_{0}(X)$ of coherent $\calox$--modules, and so the type of projective 
varieties that can carry a tilting object is severely restricted by the requirement that 
$K_{0}(X)$ be free abelian of finite rank.
\end{sit}

\begin{sit}
If the field $k$ is algebraically closed, then; see, for example, \cite[p.35f]{ARS}; any
finite-dimensional $k$--algebra $A$ is Morita--equivalent to a {\em basic\/} algebra, an algebra
$A'$ with a complete set $\{e_{i}\}_{i=1,...,N}$ of primitive orthogonal idempotents 
such that $e_{i}A'\cong e_{j}A'$ as right $A$--modules only if $i=j$.
The modules $e_{i}A'$ are then, up to $A'$--module isomorphism, the unique indecomposable 
(right) projective $A'$--modules, and, with $\rad A'$ the radical of $A'$, the modules 
$S_{i}=e_{i}A'/e_{i}\rad A'$ represent precisely the different isomorphism classes of simple 
$A'$--modules. Their respective classes form an integral basis of the Grothendieck group 
$K_{0}(A')$, isomorphic to $\ZZ^{N}$.

Further information is encoded in the {\em quiver\/} attached to $A'$, with vertices labeled
by the indices ${i=1,...,N}$, with the number of arrows from the vertex ${j}$ to the vertex ${i}$ 
equal to the (finite) dimension over $k$ of $\Ext^{1}_{A'}(S_{j},S_{i})$.

For $T$ a tilting object on $X$ and $A=\End_{X}(T)$, combining the inverse of the 
Morita--equivalence $D(A)\xto{\ \cong\ }D(A')$ with the inverse $T^{*}$ to $T_{*}$ maps 
each $e_{i}A'$ to an indecomposable direct summand $E_{i}$ of the initial tilting object $T$, 
and the direct sum $T'=\oplus_{i=1}^{N}E_{i}$ is again a tilting object on $X$. 
The difference between $T$ and $T'$ is just that $T$ may contain several copies of the same object 
$E_{i}$, that is, $T\cong \oplus_{i=1}^{N}E_{i}^{n_{i}}$, for suitable integers $n_{i}>0$. 
As the number $N$ of the pairwise nonisomorphic indecomposable summands equals the 
rank of the free abelian group $K_{0}(X)$, it is an invariant of the scheme.
\end{sit}

\begin{sit}
\label{exceptional}
In the literature, instead of the tilting object, often the set $\EE=\{E_{1},...,E_{N}\}$ of its 
indecomposable, pairwise non-isomorphic  direct summands is considered.
Many authors have studied the special case, when this set forms further 
what is also known as a {\em full, strongly exceptional collection\/} on $X$ in that in addition 
to $T=\oplus_{i}E_{i}$ being a tilting object, it is asked that the indices are (partially) ordered so that 
$\Hom_{X}(E_{j},E_{i})=0$ for $j > i$, and the summands $E_{i}$ are  furthermore required to be 
{\em simple\/}, meaning $\RR \Hom_{X}(E_{i},E_{i})\cong k[0]$ for each $i$; see, for example, 
\cite{Ru, Bon1}.
Finally, let us also mention the weaker notion of full exceptional collections,
where it is only required that $\Ext^{n}_{X}(E_{j}, E_{i})=0$ for any $n$ when $j > i$. 
In this situation; see \cite{Bon1}; the sequence $\EE$ induces on the triangulated category 
of coherent sheaves on $X$ an {\em admissible filtration\/} with layers that are semi-simple 
triangulated categories, but the whole derived category is not guaranteed to be of the form 
$D(A)$ for some algebra $A$.
\end{sit}

\begin{example}
Smooth projective varieties, at least over the complex numbers, that admit a tilting object 
include projective spaces, quadrics, Fano surfaces, various toric varieties \cite{HPe} and 
a sample of rational homogeneous varieties \cite{Boe, Sam1, Sam2}, as well as products 
of such varieties \cite{Boe}, and (iterated) projective bundles over any of these \cite{CMR}.

Dubrovin \cite[4.2.2.]{Dub} predicts in the context of complex varieties existence 
of a full, strongly exceptional collection for (smooth, projective) Fano varieties exactly when 
their {\em quantum cohomology is semi-simple\/}; see \cite{BMa} for further comments. 

For rational homogeneous manifolds $X = G/P$, with $G$ a connected complex semisimple 
Lie group, $P \subseteq G$ a parabolic subgroup, Catanese conjectures; see \cite{Boe}; that 
there should exist a tilting object, namely even a full strongly exceptional poset indexed by the 
Bruhat-Chevalley partial order of Schubert varieties in $X$.
\end{example}

\begin{remark}
If $\calt$ is equivalent to the derived category of a $k$--algebra as in \ref{algtilt}, there are 
usually many non-isomorphic, even Morita non-equivalent such algebras; see \cite{Kel}
for a more detailed discussion. However, in the situation of Theorem \ref{projtilt}, the scheme 
$X$ will often be unique up to isomorphism in view of the reconstruction theorem by Bondal 
and Orlov \cite{BO}.
\end{remark}

\subsection*{Local or Open Calabi-Yau Varieties}
\begin{sit}
\label{localCY}
Other intriguing examples, where Theorem \ref{projtilt} applies are provided by some {\em local\/} 
or {\em open Calabi-Yau varieties\/}. These include crepant resolutions of quotient 
singularities $\CC^{n}/G$, for $G$ a finite subgroup of $SL_{n}(\CC)$; see \cite[7.2.ff]{HvdB} 
for a detailed discussion of what is known or conjectured. The twisted group algebra 
$\CC[z_{1},...,z_{n}]{\ast}G$ appears then as the (suspected) endomorphism algebra of a tilting 
object. 

A second class of such examples arises as the total space of the canonical line bundle 
on those smooth projective Fano varieties that themselves carry a tilting object as foundation 
of a {\em geometric helix\/}. 
The endomorphism ring is then a ``{\em rolled-up helix algebra\/}'', a term coined by Bridgeland; 
see \cite[Thm.3.6.]{BSt} and \cite{Br} for further details.

The canonical bundle over projective space falls into both the classes just mentioned, and 
we use it as the running example to illustrate below our results. Thus, we spend a few 
lines to review this case, referring to the indicated references for details.
\end{sit}

\begin{sit}
\label{projective space}
Let $\PP=\PP^{n-1}=\PP_{k}(V)$ be the  projective space defined by an $n$--dimensional 
vector space $V$ over the field $k$, with $n\geqslant 2$. 

In Beilinson's paper \cite{Bei} that started it all, the author exhibited two tilting objects on such a projective space $\PP$, namely
\begin{align*}
T_{0} = \bigoplus_{i=0}^{n-1}\calo_{\PP}(i-n+1)\quad\text{and}\quad 
T_{1}=\bigoplus_{i=0}^{n-1}\Omega^{i}_{\PP}(i)
\end{align*}
with $\Omega^{i}_{\PP}$ the $\calo_{\PP}$-module of differential forms of degree $i$.
The associated endomorphism algebras are
\begin{align*}
A_{0} &=\End_{\PP}(T_{0}) \cong 
\bigoplus_{i,j=0}^{n-1} \Hom_{\PP}(\calo_{\PP}(j-n+1),\calo_{\PP}(i-n+1))
\cong \bigoplus_{i,j=0}^{n-1}\Sym_{i-j}(V)\\
A_{1} &=\End_{\PP}(T_{1}) \cong 
\bigoplus_{i,j=0}^{n-1} \Hom_{\PP}(\Omega^{j}_{\PP}(j),\Omega^{i}_{\PP}(i))
\cong \bigoplus_{i,j=0}^{n-1}\Lambda^{{j-i}}(V^{*})
\end{align*}
where $V^{*}$ denotes the $k$--dual vector space. Either algebra can be viewed 
as a quiver algebra on $n$ vertices labeled, say, $0,...,n-1$, with arrows from $i$ to 
$i+1$ corresponding to (a basis of) $V$, respectively $V^{*}$, and with quadratic relations 
given respectively by the kernels of the natural maps $V\otimes V\to \Sym_{2}V$ and 
$V^{*}\otimes V^{*}\to \Lambda^{2}V^{*}$. 


There are many more tilting objects, even full strongly exceptional sequences on $\PP$; see \cite{Br} for a recent discussion.
\end{sit}

\begin{sit}
Bondal \cite{Bon1, Bon2} showed that the algebras  $A_{0}$ and $A_{1}$ above are 
Koszul-duals of each other, in the sense that 
$A_{1-i}\cong A_{i}^{!} = \Ext^{\bdot}_{A_{i}}(A_{i}/\rad A_{i}, A_{i}/\rad A_{i})$ 
for $i=0,1$, where $\rad A$ denotes again the radical of the algebra $A$. Note that in 
either case, the semi-simple $k$--algebra $A_{i}/\rad A_{i}$ is just a product of $n$ copies 
of the base field $k$ with componentwise operations.
In particular, either algebra $A_{i}$ is artinian, Koszul and of finite global dimension equal to $n-1$.
\end{sit}

\begin{sit}
Let now $\pi:X\to \PP$ be the (affine) canonical bundle, the total space of the line 
bundle $\omega_{\PP}=\Omega^{n-1}_{\PP}\cong \calo_{\PP}(-n)$ over $\PP$. 
Note that this means by convention $\pi_{*}\calox \cong \Sym_{\PP}(\omega_{\PP}^{-1})$, 
thus, $X=\VV_{\PP}(\omega_{\PP}^{-1})$ in the notation of \cite{EGAII}.
The smooth variety $X$ is a {\em local\/}, also called {\em open Calabi-Yau variety\/} in that its 
canonical bundle in turn is trivial, $\omega_{X}=\Omega^{n+1}_{X}\cong \calo_{X}$.

As noted by Bridgeland \cite{Br, BSt}, any tilting object $T$ given by a full strongly exceptional 
sequence on $\PP$ pulls back to a tilting object $\pi^{*}T$ on $X$. While $\pi$ is an affine map, contracting the zero section in the affine bundle $X$ yields a projective map $p:X\to Z= \Spec K$, 
where 
\begin{align*}
K&=\oplus_{m\geqslant 0}H^{0}(\PP, (\omega_{\PP}^{-1})^{\otimes m})\cong 
\oplus_{m\geqslant 0}H^{0}(\PP, \calo_{\PP}(mn))\cong k[x_{1},...,x_{n}]^{(n)}
\end{align*}
is the $n^{th}$ Veronese subring of the polynomial ring $S=\Sym(V)\cong k[x_{1},...,x_{n}]$, 
spanned by all polynomials homogeneous of degree a multiple of $n$. If the 
characteristic of $k$ does not divide $n$ and if $k$ further contains the $n^{th}$ roots of unity, 
one may identify $K$ as well as the invariant ring under the action of the cyclic group 
$\mu_{n}$, generated by the corresponding roots of unity, acting diagonally on the variables $x_{i}$. 
That is to say $K\cong S^{\mu_{n}}$.
The ring $K$ is evidently of finite type over $k$, and so Theorem \ref{projtilt} applies.
Note that $X$ is the natural, and crepant, desingularisation of the isolated singularity of $Z$,
whence it fits also into the first class of examples mentioned in \ref{localCY}.
\end{sit}

\begin{sit}
\label{B0}
The endomorphism ring of $\pi^{*}T_{0}$ on the canonical bundle $X$ is easy to describe,
\begin{align*}
B_{0}
&=\End_{X}(\pi^{*}T_{0}) \cong \bigoplus_{i,j=0}^{n-1}S(i-j)^{(n)}
\end{align*}
that can be viewed as the algebra of $(n \times n)$--matrices, with the entries at position
$(i,j)$ sums of polynomials homogeneous of degree $i{-}j{+}mn$ for $m\geqslant 0$.

If we assume again that the characteristic of $k$ does not divide $n$ and that $k$ 
contains the corresponding roots of unity, then this algebra is isomorphic via the 
usual discrete Fourier transform to the twisted group algebra defined by the diagonal 
action of the cyclic group $\mu_{n}$ on the polynomial ring $S$, that is, 
$B_{0}\cong S{\ast}\mu_{n}$. 
This identification also exhibits $B_{0}$  clearly as a positively graded $k$--algebra, with the 
subalgebra in degree zero the semi-simple group algebra $k\mu_{n}\cong A_{0}/\rad A_{0}$.
\end{sit}

\begin{sit}
It is easily established directly that $B_{0}$ is indeed of finite global dimension equal to $n$ and homologically homogeneous, which means that all simple modules have the same projective dimension. Furthermore, it is a Calabi-Yau algebra in that it is (graded) Gorenstein with $a$--invariant 
equal to $0$ that is to say $\omega_{B_{0}}\cong B_{0}$ as graded $B_{0}$--bimodules.

The algebra $B_{0}$ is again Koszul, its Koszul-dual being the trivial extension algebra 
of $A_{1}$ by its $k$--dual $\DD(A_{1})=\Hom_{k}(A_{1},k)$, that is, 
$B_{0}^{!} \cong A_{1}\ltimes \DD(A_{1})$, an artinian symmetric algebra. In terms of representation 
by a quiver, that for $B_{0}$ is obtained from the one for $A_{0}$ by adding a copy of 
(a basis of) $V$ as additional arrows from vertex $n-1$ to vertex $0$, but keeping the same relations. 
The quiver for $B_{0}^{!}$ can be obtained from that of $A_{1}$ by here again adding a copy 
of (a basis of) $V^{*}$ as additional arrows from vertex $n-1$ to $0$ and keeping the same relations. 

As recently established by Bokland, Schedler and Wemyss  \cite{BSW}, these facts imply 
that the relations in the algebra $B_{0}$ can then be described through the derivatives of a 
single quiver (super-)potential, the simple loop in the underlying quiver of $B_{0}$ that 
corresponds to the socle element in the exterior algebra.

Bridgeland's quoted work, of which the preceding paragraph is essentially a synopsis, 
shows that the same properties are inherited by all endomorphism algebras of tilting objects on 
$X$ that come from a helix on $\PP$. In particular, the reader may want to make explicit the 
structure of $B_{1} =\End_{X}(\pi^{*}T_{1})$. 
\end{sit}

\section{Flatness and Tor-independence of Tilting Objects}
\subsection*{Tor--independence Conditions}
To simplify investigation of the behaviour of Hochschild (co-)homology under tilting, 
we will impose some flatness, or at least some $\Tor$--independence assumptions on 
tilting objects.
To this end, we make the following definition, the notion of pseudo-flatness generalised
from \cite[Defn.(80)]{BD}; see also \cite{Bu}.

\begin{defn}
Let $\calt$ be a triangulated category and assume furthermore that it is 
$R$--linear over some commutative $k$--algebra $R$.

We call a  tilting object $T$ in $\calt$  {\em flat\/} over $R$ if its endomorphism ring $A$
is flat as an $R$--algebra. The object is {\em pseudo-flat\/} if only $\Tor^{R}_{i}(A,A)=0$ 
for each $i\neq 0$.

If $\calt, \calt'$ is a pair of $R$--linear triangulated categories, $T$ a tilting 
object in $\calt$ and $T'$ a tilting object  in $\calt'$, with endomorphism 
$R$--algebras $A,A'$, respectively, then these tilting objects are {\em Tor--independent\/} 
over $R$ if $\Tor^{R}_{i}(A,A')=0$ for $i\neq 0$.
\end{defn}

Clearly, flatness implies pseudo-flatness, which in turn means that $T$ is Tor--independent
of itself, and any of these properties will hold automatically if $R$ is semi-simple, 
for example, a field.  As concerns (pseudo-)flatness, we offer the following fact that covers 
most known cases.

\begin{lemma} 
Let $f:X\to Y=\Spec R$ be a morphism from a scheme $X$ as in {\em \ref{class}\/}
 to an affine scheme $Y$ of finite type over $k$ and assume $T$ is a tilting object on $X$.
If $\calh_{i}(\calend_{X}(T)\otimes^{\LL}_{R}N)=0$ for all $i>0$ and every $R$--module $N$, 
then $T$ is a {\em flat\/}  tilting object over $R$. 

The hypthesis is satisfied in particular if $X$ is {\em flat\/} over $\Spec R$ and the 
endomorphism $\calox$--algebra $\calend_{X}(T)$ is quasiisomorphic to a 
locally free $\calox$--module, necessarily concentrated in degree $0$.
\end{lemma}

\begin{proof}
The exact functors  $\RR f_{*}(\calend_{X}(T)\otimes^{\LL}_{R}?)$ and 
$A\otimes^{\LL}_{R}?$ from $D(R)$ to itself are isomorphic by the projection formula. 
The assumption ensures that when applied to an $R$--module the first functor has 
only cohomology in non-negative degrees, while the second one always has only 
cohomology in non-positive (cohomological) degrees. Thus, $A\otimes_{R}?$ is exact 
on $R$--modules, equivalently,  $A$ is flat over $R$, whence $T$ is flat over $R$ by definition. 
\end{proof}

\begin{example}
If the tilting object $T$ is flat over $Z$  for some {\em projective\/} morphism $p:X\to Z=\Spec K$ 
to an affine scheme, then its endomorphism ring $A$ is a finite projective $K$--module, 
as it is already known to be (module-)finite over $K$ by \ref{finalg}.
\end{example}

\begin{example}
In our running example, the tilting object $\pi^{*}T_{0}$ on the anticanonical bundle $X$
over $\PP^{n}_{k}$ is not flat over the affine scheme obtained from collapsing the
zero section in $X$. However, the corresponding ring $K=S^{(n)}$ is Cohen--Macaulay,
and admits a Noether normalization, a finite morphism $\Spec K\to \Spec R$,
with $R$ smooth over $k$. One may take, for example, $R=k[x_{0}^{n},...,x_{n}^{n}]$,
the subring of the polynomial ring $S$ generated by the indicated powers of the variables.
This ring $R$ is itself a polynomial ring, and the explicit description of $B_{0}$ in \ref{B0}
shows that $B_{0}$ is a maximal Cohen-Macaulay module over $K$, thus, 
projective as an $R$--module. It follows, by definition, that $\pi^{*}T_{0}$ is flat over $R$.
\end{example}

\subsection*{Duals and Products of Tilting Objects}
Next we note that the class of tilting objects is closed under taking duals and 
``transversal products''. At least in the absolute case over an algebraically closed field,
this is certainly folklore, but we include here the details in the relative situation for
completeness.

\begin{sit}
If $X,X'$ are schemes over some common scheme $Y$, denote $p_{X}, p_{X'}$ 
the canonical projections from the fibre product $X\times_{Y}X'$ onto its factors, 
and $\LL p^{*}_{X}, \LL p^{*}_{X'}$ the respective derived inverse image functors. 
Given complexes $M,M'$ of quasi-coherent sheaves on $X,X'$ respectively, 
we set $M\boxtimes M' = \LL p_{X}^{*}M\otimes^{\LL}_{X\times_{Z}X'}\LL p_{X'}^{*}M'$.

As in Remark \ref{virtproj} we will always assume that $X,X'$ are schemes as in  
\ref{class} and that $Y$ is affine of finite type over $k$ so that the fibre product 
$X\times_{Y}X'$ is still in the category of schemes fixed in \ref{class}.
\end{sit}

\begin{proposition}
\label{dual&prod}
Let  $T$ be a  tilting object on $X$ and $T'$ a  tilting object on $X'$, respectively.  
Set $A= \End_{X}(T)$, as before, and $A'=\End_{X'}(T')$. 
\begin{enumerate}[\rm\quad(1)]
\item
\label{dualtilt}
The $\calox$--dual $T^{\vee}=\RR \calhom_{\calox}(T,\calox)$ of the perfect complex 
$T$ is again a  tilting object with $\End_{X}(T^{\vee})\cong \End_{X} (T)^{\op}$ as $K$--algebras.
In particular, $T^{\vee}$ is (pseudo-)flat along with $T$.
\item Assume $X,X'$ are flat over the affine scheme $Y=\Spec R$. 
If $T,T'$ are Tor--independent over $R$, then $T\boxtimes T'$ is a  tilting object 
on $X\times_{Y}X'$ with $\End_{X\times_{Y}X'}(T\boxtimes T')\cong A\otimes_{R}A'$.
\end{enumerate}
\end{proposition}

\begin{proof}
When restricted to the thick subcategory of perfect complexes, the functor 
$\RR\calhom_{X}(?,\calox)$ becomes an exact duality, whence $T^{\vee}$ is 
perfect and without self-extensions along with $T$. 
Now $T$ generates all of $\calt=D(X)$, so,  in particular, by Definition 
\ref{defn:tilt}(\ref{defn:tilt:genb}) the thick subcategory of perfect complexes is contained 
in $\Loc T$ and then, due to the duality and the observation in \ref{commentsgen}, that 
category is as well contained in $\Loc (T^{\vee})$. As $D(X)$ is generated by its perfect complexes, $T^{\vee}$ is also a generator. Finally, the duality $\RR\calhom_{X}(?,\calox)$ 
induces an $R$--algebra anti-isomorphism from $\End_{X}(T)$ onto $\End_{X}(T^{\vee})$, 
whence $\End_{X}(T^{\vee})\cong \End_{X} (T)^{\op}$.

As concerns (2), $T\boxtimes T'$ is a perfect complex on $X\times_{Y}X'$ as $X,X'$ are flat 
over $Y$. To confirm that this object generates, note first that due to projectivity over some 
affine scheme, $X$ carries an ample invertible sheaf, say, $\call$, and the powers 
$\call^{n}$, for $n\in\ZZ$ are contained in $\Loc T$, as that category is after all all 
of $D(X)$. The functor $?\boxtimes T'$ from $D(X)$ to 
$D(X\times_{Y}X')$ commutes with small direct sums, and so employing 
again the observation in \ref{commentsgen}, we get that $\oplus_{n\in\ZZ} \call^{n}\boxtimes T'$ 
is contained in $\Loc (T\boxtimes T')\subseteq D(X\times_{Y}X')$. 

Applying the same argument to the functor $\oplus_{n\in\ZZ} \call^{n}\boxtimes ?$
from $D(X')$ to $D(X\times_{Y}X')$ and an ample invertible sheaf 
$\call'$ on $X'$, it follows that  $\oplus_{n,m\in\ZZ} \call^{n}\boxtimes \call'^{m}$ is in turn 
contained in $\Loc \oplus_{n\in\ZZ} \call^{n}\boxtimes T'$ which we just saw to be contained 
in $\Loc T\boxtimes T'$. Now the invertible sheaf $\call\boxtimes \call'$ is ample on $X\times_{Y}X'$, whence its powers and translates generate all of $D(X\times_{Y}X')$; 
see \cite[Example 1.10]{Nee}. It follows that $T\boxtimes T'$ already generates as claimed.

It remains to verify the vanishing conditions. By the projection formula and flat base change, 
$\RHom_{X\times_{Y}X'}(T\boxtimes T', T\boxtimes T')\cong A\otimes^{\LL}_{R}A'$, 
whence the result follows from the Tor-independence of $T,T'$.
\end{proof}

\begin{remark}
The fact that $T^{\vee}$ is a tilting object along with $T$ restricts the class of algebras $A$ that occur
as endomorphism rings of tilting objects on a given scheme $X$ in that such algebras are then derived equivalent to their own opposite algebras, $D(A)\cong D(A^{\op})$, in view of Proposition \ref{dual&prod}(\ref{dualtilt}).
\end{remark}

\begin{cor}
If $X$ is flat over $Y$, and $T$ a tilting object on $X$ that is pseudo-flat over $Y$, then 
$T^{\ev} = T^{\vee}\boxtimes T$ is  tilting on $X\times_{Y}X$ with endormorphism algebra 
$A^{\ev}=A^{\op}\otimes_{R}A$.\qed
\end{cor}

To end this section, we note the following permanence property with respect to base change.
\begin{proposition}
Let $X$ be flat over the affine scheme $Y$ and $u:Y'\to Y$ an affine morphism, with $Y'$ as well
of finite type over $k$. With $u':X'=X\times_{Y}Y'\to X$ the induced morphism, if $T$ is 
a tilting object on $X$ that is flat over $Y$, then $T'=u'^{*}T$ is a tilting object on $X'$ that 
is flat over $Y'$.
\end{proposition}

\begin{proof}
Pulling back along $u'$ preserves perfection of complexes as $X$ is flat over $Y$. 
Thus, $T'$ is perfect in $D(X')$. 
The powers of any ample invertible $\calox$--module $\call$ pull back to powers of 
an ample invertible $\calo_{X'}$--module and those are contained in $u'^{*}(\Loc T)\subseteq 
\Loc u'^{*}T$. As those powers generate $D(X')$, it follows that $T'=u'^{*}T$ 
in turn generates too. Note that this is, of course, the same argument as in the proof of 
Proposition \ref{dual&prod}(2), applied to $T'\cong T\boxtimes \calo_{Y'}$.

With regard to vanishing, $p'_{*}\calend_{X'}(T') \cong p_{*}u'^{*}\calend_{X}(T) 
\cong u^{*}p_{*}\calend_{X}(T)$ by flat base change. 
Now $A=p_{*}\calend_{X}(T)$ is flat over $Y$ by assumption, whence $A'=u^{*}A$ is 
concentrated in degree zero and flat over $Y'$.
\end{proof}

\section{Hochschild Cohomology under Tilting}

\subsection*{Hochschild (Co-)Homology of Morphisms of Schemes}
If $X\to Y$ is a flat morphism of schemes or analytic spaces, then the reasonable analogue of
Hochschild theory for algebras is given by the hyper-(co-)homology of the structure sheaf 
of the diagonal in $X\times_{Y}X$. To be more precise, we recall the definition, and refer the 
reader to \cite{Swan, BFl1} for the general picture.

\begin{defn}
Let $f:X\to Y$ be a flat morphism, $\Delta:X\into X\times_{Y}X$ the embedding of the diagonal,
and denote $\calo_{\Delta}=\Delta_{*}\calox$ the structure sheaf of the diagonal.
The {\em Hochschild cohomology\/} of $X$ over $Y$ with values in a complex  $\calm$ 
from $D(X\times_{Y}X)$ is then
\begin{align*}
\Hoch^{\bdot}(X/Y,\calm) &= \Ext^{\bdot}_{X\times_{Y}X}(\calo_{\Delta},\calm)
\intertext{while the {\em Hochschild homology\/} of $X$ over $Y$ with values in $\calm$ is
defined as the hypercohomology}
\Hoch_{\bdot}(X/Y,\calm)&= H^{-\bdot}(X,\LL\Delta^{*}\calm)
\end{align*}
We write simply $\Hoch^{\bdot}(X/Y)$, respectively $\Hoch_{\bdot}(X/Y)$, if 
$\calm =\calo_{\Delta}$ is the structure sheaf of the diagonal.

In the absolute case, when $X$ is projective over the field $k$, we abbreviate further,
$\Hoch^{\bdot}(X,\calm) = \Hoch^{\bdot}(X/\Spec k, \calm)$ and 
$\Hoch^{\bdot}(X) = \Hoch^{\bdot}(X, \calo_{\Delta})$.
\end{defn}
\begin{remark}
\label{range}
If $\calm$ is a quasicoherent $\calo_{X\times_{Y}X}$--module, then the definition implies that
$\Hoch^{i}(X/Y,\calm)=0$ for $i<0$, while $\Hoch_{j}(X/Y,\calm)=0$ for $j<-\dim X$.
In general, there is no upper bound for the nonvanishing. 

However, if $X$ is locally Cohen--Macaulay and smooth over $Y$, then $\calo_{\Delta}$ is 
perfect, thus, isomorphic to a finite complex of locally free sheaves, and locally the length of such
a locally free resolution can be bounded by the relative dimension $\dim X - \dim Y$ according to the
Auslander--Buchsbaum formula. 
With $\calp$ such a locally free resolution, $\Ext^{i}_{X\times_{Y}X}(\calo_{\Delta},\calm)
\cong H^{i}(X\times_{Y}X, \calm\otimes_{X\times_{Y}X}\calp^{\vee})$, where $\calp^{\vee}$ 
denotes the $\calo_{X\times_{Y}X}$--dual of $\calp$. The hypercohomology groups on the 
right-hand side vanish for $i>2\dim X$, the sum of the length of $\calp$ and the 
dimension of $X\times_{Y}X$.
Therefore, in this case the Hochschild cohomology $\Hoch^{i}(X/Y,\calm)$ is concentrated in
the range $0\leqslant i \leqslant 2\dim X$ for any (quasi-)coherent $\calo_{X\times_{Y}X}$--module  $\calm$. Similarly, $\Hoch_{j}(X/Y,\calm)$ will be concentrated in the range 
$-\dim X\leqslant j\leqslant \dim X-\dim Y$.
\end{remark}

The following observation is somewhat pedantic but allows to exhibit clearly the action
of Hochschild cohomology on homology.
\begin{remark}
If $Y=\Spec R$ is affine, with $R$ some commutative ring, then the Hoch\-schild (co-)homology is
naturally a graded $R$--module. To make this structure explicit, note that
\begin{align*}
\Hoch_{\bdot}(X/Y,\calm)&= H^{-\bdot}(X,\LL\Delta^{*}\calm) 
\cong H^{0}(Y, H^{-\bdot}(\RR f_{*}\LL\Delta^{*}\calm))
\end{align*}
and call the complex $\RR f_{*}\LL\Delta^{*}\calm$ in $D(R)\cong D(Y)$ the 
{\em Hochschild complex\/} on $Y$ with coefficients in $\calm$.

Factoring $f:X\to Y$ through the diagonal embedding as $f:X\xto{\Delta} X\times_{Y}X\xto{f\cdot p}Y$, where $p:X\times_{Y}X\to X$ denotes any of the natural projections, 
and using the projection formula and exactness of $\Delta_{*}$, this complex can 
also be displayed as
\begin{align*}
\label{HHaction}
\tag{$\dagger$}
\RR f_{*}\LL\Delta^{*}\calm &\cong \RR(fp)_{*}\Delta_{*}\LL\Delta^{*}\calm 
\cong \RR(fp)_{*}(\calm\otimes^{\LL}_{X\times_{Y}X}\calo_{\Delta})
\end{align*}
A class $\xi\in \Hoch^{i}(X/Y)$ in Hochschild cohomology is represented by a morphism
$\xi:\calo_{\Delta}\to \calo_{\Delta}[i]$ in $D(X\times_{Y}X)$ and the induced morphism
of complexes of $R$--modules
\begin{align*}
 \RR(fp)_{*}(\calm\otimes^{\LL}_{X\times_{Y}X}\xi): \RR(fp)_{*}(\calm\otimes^{\LL}_{X\times_{Y}X}\calo_{\Delta})\to \RR(fp)_{*}(\calm\otimes^{\LL}_{X\times_{Y}X}\calo_{\Delta})[i]
\end{align*}
represents the $R$--linear action $\xi \star ?:\Hoch_{\bdot}(X/Y,\calm)\to \Hoch_{\bdot-i}(X/Y,\calm)$
of Hoch\-1529schild cohomology on homology.
\end{remark}

While the preceding definition and remark apply for any flat morphism, we now return to the situation, 
where $X$ is a scheme in the category described in \ref{class} and $Y$ is an affine scheme 
of finite type over a field $k$.

\subsection*{Preservation of the Diagonal}
With notation as in the previous section, the key result of this section is that, for a pseudo-flat 
tilting object $T$, the diagonal is preserved under tilting by $T^{\vee}\boxtimes T$, in that the 
structure sheaf of the diagonal in $X\times_{Y}X$ is transformed into $A$ with its canonical 
bimodule structure. 

To abbreviate, we set henceforth $T^{\ev}= T^{\vee}\boxtimes T$ and write accordingly the
equivalence induced by this tilting object as 
\begin{align*}
T^{\ev}_{*}= \RHom_{X\times_{Y}X}(T^{\vee}\boxtimes T,\ ):D(X\times_{Y}X)
\xto{\ \cong\  }D(A^{\ev}) 
\end{align*}

\begin{theorem}
If $X$ is flat over $Y$, and $T$ a tilting object on $X$ that is pseudo-flat over $Y$, then 
$T^{\ev}_{*}(\calo_{\Delta}) \cong A$ in $D(A^{\ev})$, where $A$ on the right is 
endowed with its canonical (right) $A^{\ev}$--module structure.
\end{theorem}

\begin{proof}
Consider the following chain of equalities and isomorphisms, where the first line simply
replaces $\calo_{\Delta}$ by its definition, and the subsequent isomorphisms result, in turn,
from the adjunction $(\LL\Delta^{*},\RR\Delta_{*}=\Delta_{*})$, Serre's ``diagonal trick''; that is, the identification of functors $\LL\Delta^{*}(-\boxtimes-)\cong -\otimes^{\LL}_{X}-$; then the adjunction $(\otimes^{\LL}_{X},
\RR\calhom_{X})$ and finally the natural identification of $T$ with $T^{\vee\vee}$:
\begin{align*}
\RHom_{X\times_{Y}X}(T^{\vee}\boxtimes T,\calo_{\Delta})
&= \RHom_{X\times_{Y}X}(T^{\vee}\boxtimes T,\Delta_{*}\calox) \\
&\cong  \RHom_{X}(\LL\Delta^{*}(T^{\vee}\boxtimes T),\calox)\\
&\cong  \RHom_{X}(T^{\vee}\otimes^{\LL}_{X} T,\calox)\\
&\cong  \RHom_{X}(T,\RR\calhom_{X}(T^{\vee},\calox))\\
&\cong  \RHom_{X}(T,T)\\
&\cong A
\end{align*}
That the identification is one of bimodules follows easily from the fact that $A$ acts (from the left)
through endomorphisms on the second factor in $T^{\vee}\boxtimes T$, while $A^{\op}$ acts (from the left) on the first one.
\end{proof}

As immediate consequences, we obtain the following results.

\begin{cor}
\label{Hochschild CoInvariance}
Assume $X$ is flat over $Y=\Spec R$, and $T$ is a pseudo-flat tilting object on it. The functor $T^{\ev}_{*}$ induces then an isomorphism of graded $R$--algebras
\begin{align}
\label{HHiso}
\tag{$*$}
\Ext^{\bdot}_{X\otimes_{Y}X}(\calo_{\Delta},\calo_{\Delta})\cong \Ext^{\bdot}_{A^{\ev}}(A,A)
\end{align}
that is, the Hochschild cohomology ring $\Hoch^{\bdot}(X/Y)$ of $X$ over $Y$ is isomorphic to the Hochschild cohomology ring $\Hoch^{\bdot}(A/R)$ of $A$ over $R$. 

Moreover, for any complex $\calm$ in $D(X\times_{Y}X)$, the same functor
induces an isomorphism 
\begin{align*}
\Ext^{\bdot}_{X\otimes_{Y}X}(\calo_{\Delta},\calm)\cong \Ext^{\bdot}_{A^{\ev}}(A,T^{\ev}_{*}\calm)
\end{align*}
of graded right modules over that ring isomorphism.
\qed
\end{cor}

\subsection*{Preservation of Smoothness}
In \cite{vdB}, van den Bergh introduced {\em smoothness\/} for an algebra over a field to mean 
that its Hochschild dimension, equal to the projective dimension of the algebra as a module 
over its enveloping algebra, is finite. Extending this definition to algebras over an arbitrary 
commutative ring, we require additionally that the algebra is {\em flat\/} over that ring. 
With this definition, we can now formulate the following improvement over 
Theorem \ref{projtilt}(\ref{fgd}).
\begin{cor}
\label{smoothness}
If $X$ and its tilting object $T$ are flat over $Y$, then $X$ is {\em smooth\/}  over $Y$ if, and 
only if, $A$ is {\em  smooth\/} over $R$.
\end{cor}

\begin{proof}
The flat morphism $X\to Y=\Spec R$ is smooth if, and only if, $\calo_{\Delta}$ is a compact 
object in $D(X\times_{Y}X)$, while $A$ is smooth over $R$ if, and only if, $A$ is flat over $R$ 
and its projective dimension over $A^{\ev}$ is finite, equivalently, $A$ is a compact object in 
$D(A^{\ev})$.
Now $A$ is flat over $R$ as $T$ is a flat tilting object, and $\calo_{\Delta}$ is compact, if, and 
only if, its image $T^{\ev}_{*}\calo_{\Delta}=A$ under the equivalence $T^{\ev}_{*}$ is compact.
\end{proof}

\begin{remark}
When algebras $A$ and $B$ of finite type over an algebraically closed field are derived equivalent, 
then finite global dimension of one implies finite global dimension of the other, but even in the 
artinian case these dimensions need not be the same. In the artinian case, finite global dimension is equivalent to smoothness as the global dimension equals the projective dimension of the algebra 
as a module over its enveloping algebra; see \cite{Hap}. We do not know whether algebras that 
appear as endomorphism rings of tilting objects on the same smooth projective scheme $X$ may 
have different global dimensions. However, using the equality of global dimension and of 
projective dimension over the enveloping algebra, we get immediately from 
Corollary \ref{Hochschild CoInvariance} that the global dimension of $A=\End_{X}(T)$ is 
bounded from below by $\max\{n| \Hoch^{n}(X)\neq 0\}$, and in all examples we are aware 
of even equality holds. 
Note that in view of Remark \ref{range} we have $\max\{i| \Hoch^{i}(X)\neq 0\}\leq 2\dim X$,
and this inequality will usually be strict.

For a concrete example, if $X=\PP^{n-1}$ as in \ref{projective space}, then 
$\max\{i| \Hoch^{i}(X)\neq 0\} = n-1 = \dim X$ equals the global dimension of either
endomorphism algebra $A_{i}$, for $i=0,1$, of the respective tilting object $T_{i}$ there.
\end{remark}

\begin{sit}
Given that the Hochschild cohomology rings of $X$ and $A$ are isomorphic, it seems 
reasonable to suspect that $X$ and $A$ have indeed as well isomorphic deformation theories. 
At least, the given tilting object lifts to any flat deformation of $X$ as it has no higher self-extensions, 
and such a lifting might conceivably still serve as a tilting object on the deformation, with 
endomorphism ring  a deformation of the original algebra $A$. However, we will not pursue this 
problem here further. Instead, we now turn to Hochschild homology.
\end{sit}

\section{Hochschild Homology and Tilting}
\subsection*{Invariance of the Hochschild complex}
\begin{theorem}
\label{invariance complex}
Assume $f:X\to Y=\Spec R$ is flat, and $T$ is a tilting object on $X$ that is 
pseudo-flat over $Y$. One then has a natural isomorphism of functors
\begin{align*}
\RR f_{*}\LL\Delta^{*}(?) \cong T^{\ev}_{*}(?)\otimes^{\LL}_{A^{\ev}}A: D(X\times_{Y}X)\lto
D(R)
\end{align*}
and for any complex $\calm$ in $D(X\times_{Y}X)$, the functor 
$T^{\ev}_{*}$ induces an isomorphism of graded $R$--modules 
$\Hoch_{\bdot}(X/Y,\calm) \cong \Hoch_{\bdot}(A/R,T^{\ev}_{*}\calm)$, linear over the isomorphism 
{\em (\ref{HHiso})\/} in Hochschild cohomology. 
In particular, $\Hoch_{\bdot}(X/Y)\cong \Hoch_{\bdot}(A/R)$.
\end{theorem}

\begin{proof}
This proof uses essentially the fact that $\RR f_{*}:D(X)\to D(R)$ admits a 
{\em right adjoint\/} $f^{!}$; see \cite[]{Nee}. Indeed, one then has the following chain of isomorphisms 
of $R$--modules, natural both in $\calm$ from $D(X\times_{Y}X)$ and $N$ from $D(R)$. The first 
one arises from the adjunctions $(\RR f_{*},f^{!})$ and $(\LL\Delta^{*},\Delta_{*})$,
\begin{align*}
\Hom_{R}(\RR f_{*}\LL\Delta^{*}\calm,N)&\cong \Hom_{X\times_{Y}X}(\calm, \Delta_{*}f^{!}N)
\intertext{the next two from applying the equivalence $T^{\ev}_{*}$ and then expanding its 
definition in the second argument,}
&\cong \Hom_{A^{\ev}}(T^{\ev}_{*}\calm, T^{\ev}_{*}\Delta_{*}f^{!}N))\\
&=  \Hom_{A^{\ev}}(T^{\ev}_{*}\calm, \RHom_{X\times_{Y}X}(T^{\vee}\boxtimes T,\Delta_{*}f^{!}N))
\intertext{while using again the adjunction $(\LL\Delta^{*},\Delta_{*})$, this time in the second 
argument, and then employing the {\em diagonal trick\/} $\LL\Delta^{*}(-\boxtimes-)
\cong -\otimes^{\LL}_{X}-$ as above, transforms this expression isomorphically into}
&\cong \Hom_{A^{\ev}}(T^{\ev}_{*}\calm, \RHom_{X}(T^{\vee}\otimes^{\LL}_{X} T, f^{!}N))
\intertext{Perfection of $T$ yields the identification $T^{\vee}\otimes^{\LL}_{X} T
\cong \RR \calhom_{X}(T,T)$, which in turn provides the isomorphism}
&\cong  \Hom_{A^{\ev}}(T^{\ev}_{*}\calm, \RHom_{X}(\RR \calhom_{X}(T,T), f^{!}N)
\intertext{Applying once again the adjunction $(\RR f_{*}, f^{!})$ together with the 
quasiisomorphisms $\RR f_{*}  \RR \calhom_{X}(T,T)\cong \RR\Hom_{X}(T,T)\cong A$ in 
$D(R)$, results in the isomorphisms}
&\cong \Hom_{A^{\ev}}(T^{\ev}_{*}\calm, \RHom_{R}(\RR f_{*}  \RR \calhom_{X}(T,T),N)\\
&\cong \Hom_{A^{\ev}}(T^{\ev}_{*}\calm, \RHom_{R}(A,N))
\intertext{and finally the adjunction $\Hom_{A^{\ev}}(U,\RHom_{R}(V,W))
\cong \Hom_{R}(U\otimes^{\LL}_{A^{\ev}}V,W)$ for complexes of (right) $A^{\ev}$--modules 
$U,V$ and (symmetric) $R$--modules $V, W$ establishes the isomorphism}
&\cong \Hom_{R}((T^{\ev}_{*}\calm)\otimes^{\LL}_{A^{\ev}}A, N)
\end{align*}
In summary, the bi-functors 
$\Hom_{R}(\RR f_{*} \LL\Delta^{*}(-),?)$ and 
$\Hom_{R}((T^{\ev}_{*}(-)\otimes^{\LL}_{A^{\ev}}A, ?)$
on $D(X\times_{Y}X)^{\op}\times D(R)$ are isomorphic, from which the first claim follows.

To see that this isomorphism of functors is linear over the isomorphism in Hochschild cohomology,
one may rewrite the above chain of isomorphisms beginning from the description
\begin{align*}
\RR f_{*}\LL\Delta^{*}\calm &\cong \RR(fp)_{*}(\calm\otimes^{\LL}_{X\times_{Y}X}\calo_{\Delta})
\end{align*}
in (\ref{HHaction}) above. The justification of each individual step in the following chain of 
isomorphisms  is the same as before, except that this time we use the right adjoint $(fp)^{!}$ 
to $\RR(fp)_{*}$,
\begin{align*}
&\Hom_{R}(\RR(fp)_{*}(\calm\otimes^{\LL}_{X\times_{Y}X}\calo_{\Delta}),N)\\
&\cong \Hom_{X\times_{Y}X}(\calm, \RR\calhom_{X\times_{Y}X}(\calo_{\Delta},(fp)^{!}N))\\
&\cong \Hom_{A^{\ev}}(T^{\ev}_{*}\calm, T^{\ev}_{*}\RR\calhom_{X\times_{Y}X}(\calo_{\Delta},(fp)^{!}N))\\
&=  \Hom_{A^{\ev}}(T^{\ev}_{*}\calm, \RHom_{X\times_{Y}X}(T^{\vee}\boxtimes T,
\RR\calhom_{X\times_{Y}X}(\calo_{\Delta},(fp)^{!}N)))\\
&\cong  \Hom_{A^{\ev}}(T^{\ev}_{*}\calm, \RHom_{R}(\RR(fp)_{*}((T^{\vee}\boxtimes T)
\otimes^{\LL}_{X\times_{Y}X}\calo_{\Delta}),N)))
\end{align*}
and, using again (\ref{HHaction}), but now on the term $ \RHom_{R}(\RR(fp)_{*}((T^{\vee}\boxtimes T)\otimes^{\LL}_{X\times_{Y}X}\calo_{\Delta})$, we find
\begin{align*}
\RR(fp)_{*}((T^{\vee}\boxtimes T)\otimes^{\LL}_{X\times_{Y}X}\calo_{\Delta})&\cong
\RR f_{*}\LL\Delta^{*}(T^{\ev})\\
&\cong T^{\ev}_{*}(T^{\ev})\otimes^{\LL}_{A^{\ev}}A\\
&\cong A^{\ev}\otimes^{\LL}_{A^{\ev}}A\\
&\cong A
\end{align*} 
whence the action of the Hochschild cohomology of $X$ through the argument $\calo_{\Delta}$ 
on the geometric side is seen to be transported into the action of the Hochschild cohomology of 
$A$ on the second factor in $(T^{\ev}_{*}\calm)\otimes^{\LL}_{A^{\ev}}A$ on the algebraic side. 
This proves that the isomorphism of functors is indeed linear over the corresponding isomorphism 
in cohomology as desired.
\end{proof}

\subsection*{Vanishing of Negative Hochschild Homology}
For algebras, Hochschild homology of modules is necessarily concentrated in non-negative
(homological) degrees, that is, $\Hoch_{i}(A/R,M)=0$ for $i<0$ and $M$ an $A^{\ev}$--module.
On the geometric side, however, as we pointed out in Remark \ref{range}, a priori, 
$\Hoch_{i}(X/Y,\calm)$ for a $\calo_{X\times_{Y}X}$--module $\calm$ is only guaranteed to 
vanish for $i < -\dim X$, due to the appearance of hypercohomology in the definition of the latter.
In the presence of a tilting object, we get better vanishing behaviour.

\begin{cor}
\label{homology vanishing}
If $f:X\to Y=\Spec R$ is flat, and $X$ admits a tilting object that is pseudo-flat over $Y$, then
\begin{enumerate}[\quad\rm (1)]
\item
\label{HH1} 
$\Hoch_{i}(X/Y)=0$ for $i<0$;
\item
\label{HH2}
$H^{i}(X,\calox)=0$ for $i\neq 0$.
\end{enumerate}
\end{cor}

\begin{proof}
Statement (\ref{HH1}) follows immediately from Theorem \ref{invariance complex} as for any 
algebra $\Hoch_{i}(A/R)\cong \Hoch_{i}(X/Y)$ vanishes for $i<0$. For (\ref{HH2}), 
note that the co-unit of the adjunction $(\LL\Delta^{*}, \Delta_{*})$ provides for a natural 
morphism $\LL\Delta^{*}(\calo_{\Delta})\cong \LL\Delta^{*}\Delta_{*}(\calox)\to \calox$ 
whose image under $\Delta_{*}$ becomes a split epimorphism as for any adjunction. 
In view of (\ref{HHaction}), this shows that $H^{i}(\RR f_{*}\calox)\cong H^{i}(X,\calox)$ 
splits off as a direct summand of $\Hoch_{-i}(X/Y)$, whence the first claim
implies the second.
\end{proof}

\subsection*{Vanishing of Hochschild Homology in the Absolute Case}
For a scheme that is projective and smooth over a field, we get an even stronger result.

\begin{theorem}
\label{vanishing absolute}
If a smooth projective scheme $X$ over a field $k$ carries a tilting object, then
$\Hoch_{i}(X)=0$ for $i\neq 0$.
\end{theorem}

\begin{proof}
Let $A=\End_{X}(T)$ be the endomorphism algebra of a tilting object on $X$ so that
$\Hoch_{\bdot}(A)\cong\Hoch_{\bdot}(X)$. As recalled before; see Theorem \ref{projtilt}(\ref{fgd}) 
and Corollary \ref{smoothness}; $A$ is then necessarily of finite global dimension, even smooth. 
A result essentially due to Keller \cite[Prop.2.5]{Kel2}; see \cite[Prop.6]{Han} for further details; 
then says that $\Hoch_{i}(A)=0$ for $i\neq 0$.
\end{proof}

\section{The Hodge Decomposition in Characteristic Zero}
We now restrict to smooth morphisms over a field of characteristic zero that we may assume to 
be the complex numbers by a suitable base change to an algebraically closed field and 
application of the Lefschetz principle. 
In this situation, Hochschild (co-)\-homology
admits a {\em Hodge decomposition\/}; see again \cite{Swan, BFl2} for details; in that
\begin{align*}
\Hoch_{i}(X/Y)&\cong \bigoplus_{p-q=i}H^{q}(X,\Omega^{p}_{X/Y})\\
\Hoch^{i}(X/Y)&\cong \bigoplus_{p+q=i}H^{q}(X, \calt^{p}_{X/Y})
\cong \bigoplus_{p+q=i}\Ext^{q}_{X}(\Omega^{p}_{X/Y}, \calox)
\end{align*}
where $\Omega^{p}_{X/Y}$ denotes the $\calox$--module of relative differential forms of 
order $p$ and $\calt^{p}_{X/Y}$ the indicated exterior power of the tangent sheaf 
$\calt_{X/Y}=\calhom_{X}(\Omega^{1}_{X/Y},\calox)$.

In the absolute case, when $Y=\Spec \CC$, Kodaira--Serre duality yields the isomorphisms
$H^{q}(X,\Omega^{p}_{X})\cong H^{n-q}(X,\Omega^{n-p}_{X})^{*}$, and so vanishing of Hochschild homology in negative degrees in Corollary \ref{homology vanishing} (\ref{HH1}) implies already
vanishing in positive degrees, reproving Theorem \ref{vanishing absolute} in this case.

\begin{theorem}
If $X$ is a smooth complex projective variety that carries a tilting object, then the Hodge 
groups satisfy $H^{q}(X,\Omega^{p}_{X})=0$ for $p\neq q$, that is, $\Hoch_{i}(X)=0$ for $i\neq 0$, 
and $\Hoch_{0}(X)\cong \oplus_{p}H^{p}(X,\Omega^{p}_{X})$.\qed
\end{theorem}

\begin{question}
How is the Hodge decomposition reflected on the algebraic side? In other words, knowing
$A=\End_{X}(T)$ for some complex projective variety $X$ with tilting object $T$, can one read off
the Hodge groups or Hodge filtration algebraically on $\Hoch_{0}(A)$?
\end{question}

In the relative situation, when $f:X\to Y$ is a smooth morphism over a field of characteristic zero,  
we get at least the following result on the Hodge modules from Corollary \ref{homology vanishing}.

\begin{theorem}
\label{Hodgevan}
If $f:X\to Y=\Spec R$ is smooth, and $X$ carries a tilting object that is pseudo-flat over $Y$, then
$H^{q}(X,\Omega^{p}_{X/Y})=0$ for $q > p$.\qed
\end{theorem}

\section{Hochschild (Co-)Homology for some Open Calabi--Yau Varieties}
\begin{sit}
Assume that the finite group $G\subset SL(V)$, for some finite-dimensional complex 
vector space $V$, {\em satisfies the crepant resolution conjecture\/} in that there is a resolution 
of singularities $X\to V/G$ with $X$ carrying a tilting object whose endomorphism ring is 
isomorphic to the twisted group ring $\calo(V){*}G$; see \cite[Ch. 7]{HvdB} 
for details and known cases and \cite{GKa} for far-reaching consequences. 
\end{sit}

\begin{sit}
The Hochschild (co-)homology of $X$ is then easily obtained, as it is known, including the
multiplicative structure on cohomology, on the algebraic side for twisted group rings; 
see e.g. \cite{Far, GKa, SWi}. 

To recall the algebraic result succinctly, we use the {\em fixed-point scheme\/} 
\[
Z =\{(v,g)\in V\times G\mid g(v) = v\}\subseteq  V\times G
\]
It is just the disjoint union $Z\cong \coprod_{g\in G}V^{g}$ of the linear spaces 
$V^{g}= \Ker(\id -g)\subseteq V$ of fixed points of the various elements $g$ in $G$. 
It is in particular naturally an affine scheme, on which the group $G$ still acts through 
$h(z,g)=(h(z), hgh^{-1})$.
The resulting quotient $Z/G$ can be identified, non canonically, with 
$\coprod_{[g]\in G{/}{\sim}}V^{g}/C_{g}$, where the disjoint union runs over the 
conjugacy classes $G{/}{\sim}$ of $G$, and $[g]\in G{/}{\sim}$ denotes the (arbitrary) 
choice of one element from each conjugacy class with $C_{g}$ its centralizer in $G$. 
This scheme is thus the disjoint union of quotient singularities resulting from the action 
of the stabilizers of the conjugacy classes in $G$ on the respective fixed linear subspace. 
It contains the original quotient singularity $V/G$ as the connected component corresponding 
to the identity element in $G$.
\end{sit}

In these terms, one has the following result, a straightforward reinterpretation of
\cite[Thm.~26 \& Thm.~31]{Far}; presented in this form in \cite{GKa} for cohomology.

\begin{theorem}
\label{HHtwistedgroup}
Given a finite-dimensional complex vector space $V$ with coordinate ring $S=\calo(V)$, 
and a finite subgroup $G\subset SL(V)$, the Hochschild (co-)homology of the twisted 
group ring $S{*}G$ is concentrated in (co-)homological degrees 
$0\leqslant i\leqslant \dim_{\CC}V$ and given, as modules over $S^{G}$, through
\begin{align*}
\Hoch_{i}(S{*}G) &\cong(\Omega^{i}_{Z})^{G}\cong 
\bigoplus_{[g]\in G/{\sim}}(\Omega^{i}_{V^{g}})^{C_{g}}
\shortintertext{and}
\Hoch^{i}(S{*}G) &\cong \Hoch_{\dim_{\CC}V-i}(S{*}G)\otimes_{\CC}(\det V)^{-1} \\
&\cong\left(\Omega^{\dim_{\CC} V - i}_{Z}\right)^{G}\otimes_{\CC}(\det V)^{-1}\\
&\cong\bigoplus_{[g]\in G/{\sim}}\left(\calt^{i-\dim_{\CC}
{V}/{V^{g}}}_{V^{g}}\otimes\det (V/V^{g})^{-1}\right)^{C_{g}}
\end{align*}
where $\calt_{V^{g}}$ denotes as before the tangent sheaf and the superscript indicates the 
appropriate exterior power of it.
\qed
\end{theorem}

On the geometric side, this then yields the following information.
\begin{cor}
If $X\to V/G$ is a crepant resolution for a finite subgroup $G\subset SL(V)$ 
so that $D(X)\simeq D(\calo(V){*}G)$, then
\begin{align*}
\Hoch_{i}(X)&\cong \bigoplus_{p}H^{p}(X,\Omega^{p+i}_{X})\cong (\Omega^{i}_{Z})^{G}
\cong \bigoplus_{[g]\in G{/}{\sim}}(\Omega^{i}_{V^{g}})^{C_{g}}\\
\shortintertext{in particular,}
\Hoch_{0}(X)&\cong \bigoplus_{p}H^{p}(X,\Omega^{p}_{X})\cong \calo(Z)^{G}\cong \prod_{[g]\in G{/}{\sim}}\calo(V^{g})^{C_{g}}
\intertext{Moreover}
\Hoch^{i}(X)&\cong \Hoch_{\dim X -i}(X)\otimes (\det V)^{-1}
\shortintertext{in particular,}
\Hoch^{0}(X)&\cong H^{0}(X,\calox)\cong \calo(V)^{G}
\end{align*}\qed
\end{cor}

\begin{question}
There remains again the question whether there is a direct interpretation of the
direct sum decomposition of Hochschild homology into the Hodge modules in terms of the 
group representation data that appear on the right-hand side, or what information relating 
those direct sum decompositions might reveal. Such interpretations 
abound for Kleinian surface singularities, the case of $G\subset SL_{2}(\CC)$, where one has 
various dictionaries provided through the McKay correspondence and its dual; see \cite{Bry} for 
a rather comprehensive geometric account. 

On the cohomological side, in general, one can interpret part of the group representation 
data as orbifold cohomology of the assumed crepant resolution if the group acts further symplectically; see \cite[Theorem 1.2.]{GKa} and the related discussion there.
\end{question}

\begin{remark}
Theorem \ref{HHtwistedgroup} extends to finite subgroups $G\subset GL(V)$, with the only modification that the factor $(\det V)^{-1}$ in the description of Hochschild cohomology has to move inside, $\Hoch^{i}(S{*}G) \cong 
\left(\Omega^{\dim_{\CC} V - i}_{Z}\otimes_{\CC}(\det V)^{-1}\right)^{G}$, and that in the 
direct sum description that follows only those conjugacy classes that lie in 
$G\cap SL(V)$ contribute; see \cite[Example 35]{Far} for details.
\end{remark}

\subsection*{The Canonical Bundle on Projective Space Revisited}
\begin{sit}
As our final example we return to the example of $X=\VV_{\PP}(\omega_{\PP}^{-1})$, the 
canonical bundle over $\PP=\PP(V)$. As before, we assume $\dim V = n\geqslant 2$, and, for simplicity, 
we work over $k=\CC$. In this situation; see \ref{B0}; the result 
$\Hoch(X) \cong \Hoch(\calo(V){*}\mu_{n})$, with $\mu_{n}$ acting diagonally through 
multiplication by roots of unity tells us in view of Theorem \ref{HHtwistedgroup} that 
\begin{align*}
\Hoch_{0}(X) &\cong \calo(V)^{\mu_{n}}\oplus\bigoplus_{1\neq g\in\mu_{n}}\CC \cong 
 \calo(V)^{\mu_{n}}\oplus \CC^{n}\\
\Hoch_{i}(X) &\cong (\Omega^{i}_{V})^{\mu_{n}}\quad\text{for $i\neq 0$}
\end{align*}
as any element $g\neq 1$ in $\mu_{n}$ leaves only the origin fixed, whence 
$Fix(g)=\{0\}$ and $\mu_{n}$ acts trivially on $\CC=\calo(Fix(g))$.
\end{sit}

\begin{sit}
Using that $\Hoch_{i}(X)=\oplus_{q\geqslant 0}H^{q}(X,\Omega^{i+q}_{X})$,
this result can easily be reconfirmed geometrically. In fact, the Zariski--Jacobi 
sequence for the smooth affine structural morphism $\pi:X\to \PP$ is
\begin{align*}
0\to \pi^{*}\Omega^{1}_{\PP}\to \Omega^{1}_{X}\to \Omega^{1}_{X/\PP}\to 0
\end{align*}
and, as $\pi_{*}\calox = \Sym_{\PP}(\omega_{\PP}^{-1})$, the $\calox$--module of relative 
differentials $\Omega^{1}_{X/\PP}$ can be identified with $\pi^{*}(\omega_{\PP}^{-1})$. 
Taking exterior powers yields for any $p\geqslant 0$ a short exact sequence of $\calox$--modules 
 \begin{align*}
0\to \pi^{*}\Omega^{p}_{\PP}\to \Omega^{p}_{X}\to 
\pi^{*}(\Omega^{p-1}\otimes_{\PP}\omega_{\PP}^{-1})\to 0
\end{align*}
Taking into account that $H^{\bdot}(X,\ ) \cong H^{\bdot}(\PP, \pi_{*}(\ ))$ and that 
$H^{q}(X,\Omega^{p}_{X})=0$ for $q>p$ by Theorem \ref{Hodgevan}, the resulting long 
exact cohomology sequence becomes
\[
\xymatrix@R-=2pt@C-2pt{
0\ar[r]&H^{0}(\PP,\pi_{*}\pi^{*}\Omega^{p}_{\PP})\ar[r]&H^{0}(X,\Omega^{p}_{X})\ar[r]&
H^{0}(\PP,\pi^{*}\pi_{*}\Omega^{p-1}_{\PP}\otimes\omega_{\PP}^{-1})\ar[r]&\\
\cdots\ar[r]&H^{i}(\PP,\pi_{*}\pi^{*}\Omega^{p}_{\PP})\ar[r]&H^{i}(X,\Omega^{p}_{X})\ar[r]&
H^{i}(\PP,\pi^{*}\pi_{*}\Omega^{p-1}_{\PP}\otimes\omega_{\PP}^{-1})\ar[r]&\\
\cdots\ar[r]&H^{p+1}(\PP,\pi_{*}\pi^{*}\Omega^{p}_{\PP})\ar[r]&0
}
\]
where the occurring tensor products are taken over $\calo_{\PP}$. The first three terms 
form a short exact sequence, isomorphic to 
\[
0\to\bigoplus_{m\geqslant 0}H^{0}(\PP,\Omega^{p}_{\PP}(mn))\to
H^{0}(X,\Omega^{p}_{X})\to
\bigoplus_{m\geqslant 1}H^{0}(\PP,\Omega^{p-1}_{\PP}(mn))\to 0
\]
It can be identified with the corresponding short exact sequence that results from 
restricting the Koszul complex over $V$ to degrees that are multiples of $m$,
\[
0\to\bigoplus_{m\geqslant 0}H^{0}(\PP,\Omega^{p}_{\PP}(mn))\to
(\Omega^{p}_{V})^{(n)}\to
\bigoplus_{m\geqslant 1}H^{0}(\PP,\Omega^{p-1}_{\PP}(mn))\to 0
\]
and thereby yields the isomorphism
\begin{align*}
(\Omega^{p}_{V})^{(n)}\cong H^{0}(X,\Omega^{p}_{X})\subseteq \Hoch_{p}(X)
\end{align*}
from the invariant differential forms on $V$ onto the indicated direct summand of the
Hochschild homology of $X$. Using that for $q\geqslant 1, m\geqslant 0,$ the cohomology groups
$H^{q}(\PP,\Omega^{p}_{\PP}(m))$ vanish except for $p=q$ and $m=0$, in which case
$H^{p}(\PP,\Omega^{p}_{\PP})=\CC$, we reaffirm the result stated above. In particular, 
the $\calo(V)^{G}$--torsion submodule of $\Hoch_{0}(X)$ appears once as
$\oplus_{1\neq g\in\mu_{n}}\CC$ and second as $\oplus_{p=1}^{n-1}H^{p}(X,\Omega^{p}_{X})$.
\end{sit}

\end{document}